\documentclass[11pt]{article}
\usepackage[usenames]{color}
\usepackage{amssymb}
\usepackage[colorlinks=true,
linkcolor=webgreen,
filecolor=webbrown,
citecolor=webgreen]{hyperref}

\definecolor{webgreen}{rgb}{0,.5,0}
\definecolor{webbrown}{rgb}{.6,0,0}

\newtheorem{theorem}{Theorem}
\newtheorem{lemma}{Lemma}

\def\slfrac#1#2{\hbox{\kern.1em %
 \raise.5ex\hbox{\the\scriptfont0 #1}\kern-.11em %
 /\kern-.15em\lower.25ex\hbox{\the\scriptfont0 #2}}}

\newcommand{\bsq}{{\vrule height .9ex width .8ex depth -.1ex }}
\newcommand\red[1]{\mbox{\color{red} #1}}
\newcommand{\eeq}{\end{equation}}
\newcommand{\beql}[1]{\begin{equation}\label{#1}}
\newcommand{\eqn}[1]{(\ref{#1})}

\newcommand{\ba}{{\mathbf {a}}}
\newcommand{\bs}{{\mathbf {s}}}
\newcommand{\bc}{{\mathbf {c}}}

\newcommand{\faf}{\mbox{\boldmath $\alpha$}}
\newcommand{\fbe}{\mbox{\boldmath $\beta$}}
\newcommand{\NN}{{\mathbb {N}}}
\newcommand{\PP}{{\mathbb {P}}}

\makeatletter
\def\@sect#1#2#3#4#5#6[#7]#8{\ifnum #2>\c@secnumdepth
     \def\@svsec{}\else
     \refstepcounter{#1}\edef\@svsec{\csname the#1\endcsname.\hskip .75em }\fi
     \@tempskipa #5\relax
      \ifdim \@tempskipa>\z@
        \begingroup #6\relax
          \@hangfrom{\hskip #3\relax\@svsec}{\interlinepenalty \@M #8\par}%
        \endgroup
       \csname #1mark\endcsname{#7}\addcontentsline
         {toc}{#1}{\ifnum #2>\c@secnumdepth \else
                      \protect\numberline{\csname the#1\endcsname}\fi
                    #7}\else
        \def\@svsechd{#6\hskip #3\@svsec #8\csname #1mark\endcsname
                      {#7}\addcontentsline
                           {toc}{#1}{\ifnum #2>\c@secnumdepth \else
                             \protect\numberline{\csname the#1\endcsname}\fi
                       #7}}\fi
     \@xsect{#5}}
\def\@begintheorem#1#2{\it \trivlist \item[\hskip \labelsep{\bf #1\ #2.}]}
\makeatother

\begin{document}
\begin{center}
{\Large{\bf Numerical Analogues of Aronson's Sequence}} \\
\vspace{1\baselineskip}
{\em Benoit Cloitre} \\
13 rue Pinaigrier \\
Tours 37000, FRANCE \\
(Email: abcloitre@wanadoo.fr) \\
\vspace{1\baselineskip}
{\em N. J. A. Sloane} \\
AT\&T Shannon Labs \\
 Florham Park, NJ 07932--0971, USA \\
(Email: njas@research.att.com) \\
\vspace{1\baselineskip}
{\em Matthew J. Vandermast} \\
53 Piaget Avenue \\
Clifton, NJ 07011--1216, USA \\
(Email: ghodges14@msn.com) \\
\vspace{2\baselineskip}
March 28, 2003 \\
\vspace{2\baselineskip}
{\bf Abstract} \\
\vspace{.5\baselineskip}
\end{center}

\setlength{\baselineskip}{1.5\baselineskip}

Aronson's sequence 1, 4, 11, 16, $\ldots$ is defined by the English sentence
``t is the first, fourth, eleventh, sixteenth, $\ldots$ letter
of this sentence.''
This paper introduces some numerical analogues, such as:
$a(n)$ is taken to be the smallest positive integer greater
than $a(n-1)$ which is consistent with the condition
``$n$ is a member of the sequence if and only if $a(n)$ is odd.''
This sequence can also be characterized by its ``square'',
the sequence $a^{(2)} (n) = a(a(n))$, which equals $2n+3$ for $n \ge 1$.
There are many generalizations of this sequence,
some of which are new, while others throw new light
on previously known sequences.

\section{Introduction}
Aronson's sequence, given in the Abstract, 
is a classic example of a self-referential sequence
(\cite{Aro85}, \cite{Hof85},
sequence M3406 in \cite{EIS}, 
\htmladdnormallink{A5224}{http://www.research.att.com/cgi-bin/access.cgi/as/njas/sequences/eisA.cgi?Anum=A005224} in \cite{OEIS}).
It is somewhat unsatisfactory because of the ambiguity in
the English names for numbers over 100 --- some people say
``one hundred and one'', while
others say ``one hundred one.''
Another well-known example is Golomb's sequence, 
in which the
$n^{\rm th}$
term $G(n)$ (for $n \ge 1$) is the number of times $n$ appears
in the sequence (\htmladdnormallink{A1462}{http://www.research.att.com/cgi-bin/access.cgi/as/njas/sequences/eisA.cgi?Anum=A001462} in \cite{OEIS}):
$$1,2,2,3,3,4,4,4,5,5,5,6,6,6,6,7,7,7,7,8, \ldots$$
There is a simple formula for $G(n)$: it is
the nearest integer to (and approaches)
$$\phi^{2- \phi} n^{\phi -1} ~,$$
where $\phi = (1+ \sqrt{5})/2$ (\cite{Gol66},
\cite[Section E25]{UPNT}).

Additional examples can be found in Hofstadter's books \cite{GEB},
\cite{Hof85} and in \cite{UPNT} and \cite{OEIS}.
However, the sequence $\{a(n)\}$ mentioned in the Abstract appears to be new,
as do many of the other sequences we will discuss.
We will also give new properties of some sequences that
have been studied elsewhere.

Section~\ref{Sec2} discusses the sequence mentioned in the Abstract,
and also introduces the ``square'' of a sequence.
Some simple generalizations (non-monotonic,
``even'' and ``lying'' versions) are described in
Section~\ref{Sec3}.
The original sequence is
based on examination of the sequence modulo 2.
In Section~\ref{Sec4} we consider various ``mod $y$''
generalizations.
Section~\ref{Sec5} extends both the original sequence
and the ``mod $y$'' generalizations by defining
the ``Aronson transform'' of a sequence.
Finally, Section~\ref{Sec6} briefly considers
the case when the rule defining the sequence
depends on more than one term.

There are in fact a large number of possible generalizations
and we shall only mention some of them.
We have not even analyzed all the sequences that we do mention.
In some cases we just list the first few terms and
invite the reader to investigate them himself.
We give the identification numbers of
these sequences in \cite{OEIS} --- the entries there will be 
updated as more information becomes available.

We have also investigated sequences arising when \eqn{Eq2}
is replaced by the following rule: $s(1) = 1, s(n) = s(n-1) + y$
if $n$ is already in the
sequence, $s(n) = s(n-1) + z$ otherwise, for specified
values of $x, y, z$. This work
will be described elsewhere \cite{CSV03}.

\paragraph{Notation.} ``Sequence'' here usually means an infinite sequence
of nonnegative numbers. ``Monotonically increasing''
means that each term is strictly greater than the previous term.
$\PP = \{1,2,3, \ldots \}$,
$\NN = \{0,1,2,3, \ldots \}$.
 
\section{$n$ is in sequence if and only if $a(n)$ is odd}\label{Sec2}
Let the sequence $a(1)$, $a(2)$, $a(3), \ldots$ be defined
by the rule that $a(n)$ is the smallest positive integer
$>$ $a(n-1)$ which is consistent with the condition that 
\beql{Eqa}
\mbox{``$n$ is a member of the sequence if and only if $a(n)$ is odd.''}
\eeq

The first term, $a(1)$, could be 1, since 1 is odd
and 1 would be in the sequence.
It could also be 2, since then 1 would not be in the
sequence (because the terms must increase) and 2 is even.
But we must take the {\em smallest} possible value, so $a(1) =1$.
Now $a(2)$ cannot be 2, because 2 is even.
Nor can $a(2)$ be 3, for then 2 would not be in the sequence but $a(2)$
would be odd.
However, $a(2) =4$ is permissible, so we {\em must}
take $a(2) =4$, and then 2 and 3 are not in the sequence.

So $a(3)$ must be even and $>4$, and $a(3) =6$ works.
Now 4 is in the sequence, so $a(4)$ must be odd, and $a(4) =7$ works.
Continuing in this way we find that the first few terms are as follows
(this is \htmladdnormallink{A79000}{http://www.research.att.com/cgi-bin/access.cgi/as/njas/sequences/eisA.cgi?Anum=A079000}):
$$
\begin{array}{rrrrrrrrrrrrrr}
n: & 1 & 2 & 3 & 4 & 5 & 6 & 7 & 8 & 9 & 10 & 11 & 12 & \cdots \\
a(n): & 1 & 4 & 6 & 7 & 8 & 9 & 11 & 13 & 15 & 16 & 17 & 18 & \cdots
\end{array}
$$

Once we are past $a(2)$ there are no further complications,
$a(n-1)$ is greater than $n$, and we {\em can},
and therefore {\em must}, take
\beql{Eq2}
a(n) = a(n-1) + \epsilon ~,
\eeq
where $\epsilon$ is 1 or 2 and is given by:
$$
\begin{array}{ccc}
~ & \mbox{$a(n-1)$ even} & \mbox{$a(n-1)$ odd} \\ [+.1in]
\mbox{$n$ in sequence} & 1 & 2 \\
\mbox{$n$ not in sequence} & 2 & 1
\end{array}
$$
The gap between successive terms for $n \ge 3$ is either 1 or 2.

The analogy with Aronson's sequence is clear.
Just as Aronson's sentence indicates exactly which of its terms are t's,
$\{a(n)\}$ indicates exactly which of its terms are odd.

We proceed to analyze the behavior of this sequence.

First, all odd numbers $\ge 7$ occur.
For suppose $2t+1$ is missing.
Therefore $a(i) = 2t$, $a(i+1) = 2t+2$ for some $i \ge 3$.  From
the definition, this means $i$ and $i+1$ are missing,
implying a gap of at least 3, a contradiction.

Table \ref{T1} shows the first 72 terms, with the even numbers colored red.
\begin{table}[htb]
\caption{The first 72 terms of the sequence ``$n$ is in sequence if and only if $a(n)$ is odd.''}
$$
\begin{array}{rrrrrrrrrrrr}
n: & 1 & \multicolumn{1}{r|}{2} & 3 & 4 & 5 & 6 & 7 & \multicolumn{1}{r|}{8} & 9 & 10 \\
a(n): & 1 & \multicolumn{1}{r|}{\red{4}} & \red{6} & 7 & \red{8} & 9 & 11 & \multicolumn{1}{r|}{13} & 15 & \red{16} \\
~ & ~ \\
n: & 11 & 12 & 13 & 14 & 15 & 16 & 17 & 18 & 19 & \multicolumn{1}{r|}{20} \\
a(n): & 17 & \red{18} & 19 & \red{20} & 21 & 23 & 25 & 27 & 29 & \multicolumn{1}{r|}{31} \\
~ & ~ \\
n: & 21 & 22 & 23 & 24 & 25 & 26 & 27 & 28 & 29 & 30 \\
a(n): & 33 & \red{34} & 35 & \red{36} & 37 & \red{38} & 39 & \red{40} & 41 & \red{42} \\
~ & ~ \\
n: & 31 & 32 & 33 & 34 & 35 & 36 & 37 & 38 & 39 & 40 \\
a(n): & 43 & \red{44} & 45 & 47 & 49 & 51 & 53 & 55 & 57 & 59 \\
~ & ~ \\
n: & 41 & 42 & 43 & \multicolumn{1}{r|}{44} & 45 & 46 & 47 & 48 & 49 & 50 \\
a(n): & 61 & 63 & 65 & \multicolumn{1}{r|}{67} & 69 & \red{70} & 71 & \red{72} & 73 & \red{74} \\
~ & ~ \\
n: & 51 & 52 & 53 & 54 & 55 & 56 & 57 & 58 & 59 & 60 \\
a(n): & 75 & \red{76} & 77 & \red{78} & 79 & \red{80} & 81 & \red{82} & 83 & \red{84} \\
~ & ~ \\
n: & 61 & 62 & 63 & 64 &  65 & 66 & 67 & 68 & 69 & 70 \\
a(n): & 85 & \red{86} & 87 & \red{88} & 89 & \red{90} & 91 & \red{92} & 93 & 95 \\
~ & ~ \\
n: & 71 & 72 & \multicolumn{2}{l}{\cdots} \\
a(n): & 97 & 99 & \multicolumn{2}{l}{\cdots} \\
\end{array}
$$
\label{T1}
\end{table}

Examining the table, we see that there are three consecutive numbers,
6, 7, 8, which are necessarily followed by three consecutive odd numbers,
$a(6)=9$, $a(7)=11$, $a(8)=13$.
Thus 9 is present, 10 is missing, 11 is present, 12 is missing,
and 13 is present.
Therefore the sequence continues with $a(9) =15$ (odd),
$a(10) = 16$ (even), $\ldots$, $a(13)=19$ (odd), $a(14) = 20$ (even).
This behavior is repeated for ever.
A run of consecutive numbers is immediately followed by
a run of the same length of consecutive odd numbers.

Let us define the $k^{\rm th}$ segment (for $k \ge 0$)
to consist of the terms $a(n)$ with $n= 9 \cdot 2^k - 3 + j$
where $-3 \cdot 2^k \le j \le 3 \cdot 2^k -1$.
In the table the segments are separated by vertical lines.
The first half of each segment, the terms where $j \le 0$,
consists of consecutive numbers given by $a(n) = 12 \cdot 3^k -3 + j$;
the second half,
where $j \ge 0$, consists of consecutive odd numbers given by
$a(n) = 12 \cdot 3^k -3 + 2j$.
We can combine these formulae,
obtaining an explicit description for the sequence:
$$
a(1) = 1, \quad a(2) =4 ,
$$
and subsequent terms are given by
\beql{Eq1}
a(9 \cdot 2^k - 3+j) = 12 \cdot 2^k -3 + \frac{3}{2} j + \frac{1}{2} |j|
\eeq
for $k \ge 0$, $-3\cdot 2^k \le j < 3\cdot 2^k$.

The structure of this sequence is further revealed by
examining the sequence of first differences,
$\Delta a(n) = a(n+1) - a(n)$,
$n \ge 1$, which is
\beql{Eq1d}
3,2,1,1,1,2,2,2,1^6, 2^6,1^{12}, 2^{12}, 1^{24}, 2^{24}, \ldots
\eeq
(\htmladdnormallink{A79948}{http://www.research.att.com/cgi-bin/access.cgi/as/njas/sequences/eisA.cgi?Anum=A079948}), where we have written $1^m$ to indicate a string of $m$ 1's, etc.
The oscillations double in length at each step.

Segment 0 begins with an even number,
6, but all other segments begin with an odd number,
$9 \cdot 2^k -3$.
All odd numbers occur in the sequence except 3 and 5.
The even numbers that occur are 4, 6, 8 and all numbers $2m$ with
$$9 \cdot 2^{k-1} -1 \le m \le 6 \cdot 2^k -2 , \quad k \ge 1 \,.$$

The sequence of differences, \eqn{Eq1d},
can be constructed from the words in a certain formal
language (cf. \cite{Lot83}).
The alphabet is $\{1, 2, 3\}$ and we define morphisms
$\theta(1) = 2, \theta(2) = 1,1$. Then \eqn{Eq1d} is the concatenation 
\beql{Eq1s}
S_{-1}, S_0, S_1, S_2, \ldots ~ ,
\eeq
where
\beql{Eq1t}
S_{-1} = \{3,2\}, S_0 = \{1,1,1\},
S_{k+1} = \theta (S_k) \mbox{~for~} k \ge 0 ~.
\eeq
To prove this, note that for $n \ge 3$, a difference of $2$ only
occurs in $\{a(n)\}$ between a pair of odd numbers.
Suppose $a(i) = 2j+1, a(i+1) = 2j+3$; then $a(2j+1)
= 2x+1$ (say), $a(2j+2) = 2x+2, a(2j+3) = 2x+3$,
producing two differences of $1$. Similarly, if there is a difference of $1$,
say $a(i)=j, a(i+1)=j+1$, then $a(j)=2x+1, a(j+1)=2x+3$,
a difference of $2$.

The ratio $n/a(n)$, which is the fraction of numbers 
that are in the sequence,
rises from close to $2/3$ at the beginning of segment 
$k$ (assuming $k$ is large), reaches a maximum $3/4$ at
the midpoint of the segment, then falls back to $2/3$ at 
the end of the segment.
It is not difficult to show that if $n$ is chosen at
random in the
$k^{\rm th}$
segment then the average value
of the fraction of numbers in the sequence at that point approaches
$$\frac{3}{4} - \frac{1}{4} \log \frac{32}{27} = 0.7075 \ldots$$
for large $k$.

The sequence has an alternative characterization in terms of its ``square.''

The {\em square} of a sequence $\bs = \{s(n) : n \ge n_0 \}$ is given by
$\bs^{(2)} = \{s(s(n)) : n \ge n_0 \}$.
If $\bs$ is monotonically increasing so is $\bs^{(2)}$.

\begin{lemma}\label{L1}
Let $\bs$ be monotonically increasing.
Then $n~ (\ge n_0)$ is in the sequence $\bs$ if
and only if $s(n)$ is in the sequence $\bs^{(2)}$.
\end{lemma}

\paragraph{Proof.}
If $n$ is in the sequence, $n=s(i)$ for some $i \ge n_0$,
and $s(n) = s(s(i))$ is in $\bs^{(2)}$.
Conversely, if $s(n) \in \bs^{(2)}$, $s(n) = s(s(i))$ for
some $i \ge n_0$, and since $\bs$ is monotonically increasing,
$n = s(i)$.~~~$\bsq$

For our sequence $\ba = \{a(n)\}$, examination of Table \ref{T1} shows that
$\ba^{(2)} = \{1,5,7,9,11, \ldots \} = \{1\} \cup 2 \PP +3$.
This can be used to characterize $\ba$.
More precisely, the sequence can be defined by:
$a(1) =1$, $a(2) =4$, $a(3) =6$ and, for $n \ge 4$,
$a(n)$ is the smallest positive integer which is consistent 
with the sequence being monotonically increasing and satisfying 
$a(a(n)) = 2n+3$ for $n \ge 2$.

This is easily checked.
Once the first three terms are specified, the 
rule $a(a(n)) = 2n+3$ determines the remaining terms uniquely.

In fact that rule also forces $a(2)$ to be 4,
but it does not determine $a(3)$, since there is an earlier
sequence $\{a'(n)\}$ (in the lexicographic sense) satisfying $a' (1) =1$,
$a' (a'(n)) = 2n+3$ for $n \ge 2$, namely
$$1,4,5,7,9,10,11,12,13,15,17,19,21,22, \ldots,$$
(\htmladdnormallink{A80596}{http://www.research.att.com/cgi-bin/access.cgi/as/njas/sequences/eisA.cgi?Anum=A080596}), and given by $a' (1) =1$,
\beql{Eq20}
a' (6 \cdot 2^k - 3 +j) = 8 \cdot 2^k -3 + \frac{3}{2} j + \frac{1}{2} |j|
\eeq
for $k \ge 0$, $-2^{k+1} \le j < 2^{k+1}$.

As the above examples show, the square of
a sequence does not in general determine the sequence uniquely.
A better way to do this is provided by the ``inverse Aronson transform'',
discussed in Section~\ref{Sec5}.

\section{First generalizations}\label{Sec3}
The properties of $\{a(n)\}$ given in Section \ref{Sec2}
suggest many generalizations, some of which
will be discussed in this and the following sections.

\paragraph{(3.1) Non-monotonic version.}
If we replace ``$a(n) > a(n-1)$'' in the 
definition by ``$a(n)$ is not already in the 
sequence'', we obtain a completely different sequence, suggested 
by J.~C. Lagarias (personal communication):
$b(n)$, $n \ge 1$,
is the smallest positive integer not already in 
the sequence which is consistent with the condition 
that ``$n$ is a member of the sequence if and only if $b(n)$ is odd.''
This sequence (\htmladdnormallink{A79313}{http://www.research.att.com/cgi-bin/access.cgi/as/njas/sequences/eisA.cgi?Anum=A079313}) begins:
$$
\begin{array}{rrrrrrrrrrr}
n: & 1 & 2 & 3 & 4 & 5 & 6 & 7 & 8 & 9 & 10 \\
b(n): & 1 & 3 & 5 & \red{2} & 7 & \red{8} & 9 & 11 & 13 & \red{12} \\
~ & ~ \\ [+.1in]
n: & 11 & 12 & 13 & 14 & 15 & 16 & 17 & 18 & 19 & 20 \\
b(n): & 15 & 17 & 19 & \red{16} & 21 & 23 & 25 & \red{20} & 27 & 29
\end{array}
$$
The even
members are shown in red.
The behavior is simpler than that of $\{a(n)\}$,
and we leave it to the reader to show that, for $n \ge 5$,
$b(n)$ is given by
\begin{eqnarray*}
b(4t-2) & = & 4t \,, \\
b(4t-1) & = & 6t -3 \,, \\
b(4t) & = & 6t-1 \,, \\
b(4t+1) & = & 6t +1 \,.
\end{eqnarray*}
All odd numbers occur.
The only even numbers are 2 and $4t$, $t \ge 2$.
(The square ${\mathbf {b}}^{(2)}$ is not so interesting.)

\paragraph{(3.2) ``Even'' version.}
If instead we change ``odd'' in the definition of $\{a(n)\}$ to ``even'',
we obtain a sequence $\bc$ which is best started at $n=0$:
$c(n)$, $n \ge 0$, is the smallest nonnegative integer $> c(n-1)$ 
which is consistent with the condition that
\beql{Eqc}
\mbox{``$n$ is a member of the sequence if and only if $c(n)$ is even.''}
\eeq
This is \htmladdnormallink{A79253}{http://www.research.att.com/cgi-bin/access.cgi/as/njas/sequences/eisA.cgi?Anum=A079253}:
0, 3, 5, 6, 7, 8, 10, 12, 14, 15, $\ldots$.
It is easily seen that $c(n) = a(n+1) -1$ for $n \ge 0$,
so there is nothing essentially new here.
Also $\bc^{(2)} = \{0\} \cup 2 \PP +4$.

\paragraph{(3.3) The ``lying'' version.}
The lying version of Aronson's sequence is based on
the completely false sentence ``t is the 
second, third, fifth, sixth, 
seventh, $\ldots$ letter of this sentence.''
The sentence specifies exactly those letters that are 
not t's, and produces the sequence (\htmladdnormallink{A81023}{http://www.research.att.com/cgi-bin/access.cgi/as/njas/sequences/eisA.cgi?Anum=A081023}) 2, 
3, 5, 6, 7, 8, 9, 10, 11, 12, $\ldots$.

Just as $\{a(n)\}$ is an analogue of the original sequence, we can define
an analogue
$\{d(n): n \ge 1 \}$ of this sequence 
by saying that: $d(n)$ is the smallest positive 
integer $> d(n-1)$ such that the condition ``$n$ is in the sequence
if and only if $d(n)$ is odd'' is false.
Equivalently, the condition ``either $n$ is in the sequence and
$d(n)$ is even or $n$ is not in the sequence and $d(n)$ is odd''
should be true.
The resulting sequence
(\htmladdnormallink{A80653}{http://www.research.att.com/cgi-bin/access.cgi/as/njas/sequences/eisA.cgi?Anum=A080653}) begins 2, 4, 5, 6, 8, 10, 11, 12, 13, 14, $\ldots$.
We will give an explicit formula for $d(n)$ in
the next section.

A related sequence is also of interest. Let $\{d'(n)\}$
be defined by $d'(1) = 2$, and, for $n>1$, 
$d'(n)$ is the smallest integer greater than $d'(n-1)$
such that the condition ``$n$ and $d'(d'(n))$ have opposite parities''
can always be satisfied.
One can show that this is the sequence
$$
2, 4, 5, 7, 8, 9, 11, 12, 13, 14, 16, \ldots
$$
(\htmladdnormallink{A14132}{http://www.research.att.com/cgi-bin/access.cgi/as/njas/sequences/eisA.cgi?Anum=A014132}), the complement of the triangular numbers,
with $d'(n) ~=~ n ~+~ \mbox{nearest integer to~} \sqrt{2n}$.

\section{The ``mod $m$'' versions}\label{Sec4}
Both $\{a(n)\}$ and $\{c(n)\}$ are defined modulo 2.
Another family of generalizations is based on replacing 
2 by some fixed integer $y \ge 2$.
To this end we define a sequence $\{s(n): 
n \ge n_0 \}$ by specifying a starting 
value $s(n_0) = s_0$, and the condition that 
$n$ is in the sequence if and only 
if $s(n) \equiv z$ $(\bmod~y)$, where $y$ and $z$ are given.

Although we will not digress to consider this 
here, it is also of interest to see 
what happens when ``if and only if'' in 
the definition is replaced by either ``if'' or ``only if.''
(We mention just one example. The 
above sequence $\{d(n)\}$, prefixed by $d(0) = 0$, 
can be defined as follows: $d(n)$ is the smallest nonnegative number
$> d(n-1)$ such that the condition ``$n~ (n \ge 0)$ is in the sequence
only if $d(n)$ is even'' is satisfied.)

We saw in the previous section that $\{a(n)\}$ 
can also be characterized by the property that 
its square $a^{(2)} (n) = a(a(n))$ is equal 
to $2n+3$ for $n \ge 2$ (together with 
some appropriate initial conditions).
This too can be generalized by specifying that 
the sequence $\{s(n)\}$ satisfy $s(s(n)) = yn+z$, for 
given values of $y$ and $z$.
The two generalizations are related, but usually lead 
to different sequences.
The $s(s(n))$ family of generalizations will connect the 
present investigation with several sequences that have already 
appeared in the literature.
There are too many possibilities for us to 
give a complete catalogue of all the sequences 
that can be obtained from these generalizations.
Instead we will give a few key examples and one general theorem.
Many other examples can be found in \cite{OEIS}.

A simple ``mod~3'' generalization is: $e(1)=2$, and, for 
$n > 1$, $e(n)$ is the smallest integer 
$>e(n-1)$ which is consistent with the condition that
\beql{Eqe}
\begin{array}{l}
\mbox{``$n$ is a member of the sequence if and only if
$e(n)$ is a multiple of 3.''}
\end{array}
\eeq
This turns out to be James Propp's sequence
$$2,3,6,7,8,9,12,15,18,19, \ldots ~,$$
which appeared as sequence M0747 in \cite{EIS}
(\htmladdnormallink{A3605}{http://www.research.att.com/cgi-bin/access.cgi/as/njas/sequences/eisA.cgi?Anum=A003605}in \cite{OEIS}).
Propp gave a different (although equivalent) definition
involving the square of the sequence:
$\{e(n)\}$ is the unique monotonically increasing sequence satisfying
$e(e(n)) =3n$ for all $n \ge 1$.
Michael Somos (personal communication) observed that this sequence satisfies
\begin{eqnarray}\label{Eq46}
e(3n) & = & 3e(n) , \nonumber \\
e(3n+1) & = & 2e(n) + e(n+1) , \nonumber \\
e(3n+2) & = & e(n) + 2e(n+1) .
\end{eqnarray}
An analysis similar to that for $\{a(n)\}$ leads to the following
explicit formula, which appears to be new:
\beql{Eq47}
e(2 \cdot 3^k +j) = 3^{k+1} + 2j + |j| \,,
\eeq
for $k \ge 0$ and $-3^k \le j < 3^k$.

A sequence closely related to $\{e(n)\}$ had earlier been studied by
Arkin et~al. \cite[Eq. (12)]{Ark90}.
This is the sequence $e' (n) = e(n) -n$ (\htmladdnormallink{A6166}{http://www.research.att.com/cgi-bin/access.cgi/as/njas/sequences/eisA.cgi?Anum=A006166}).
Arkin et~al. give a recurrence similar to \eqn{Eq46}.

The sequence $\{e(n)\}$ can be generalized as follows.

\begin{theorem}\label{Th1}
Let $y$ and $z$ be integers of opposite parity satisfying
\beql{Eq48a}
y \ge 2 ,~ y +z \ge 1,~ 2y +z \ge 4 \,.
\eeq
Then there is a unique monotonically increasing sequence $\{f(n)\}$ 
satisfying
$f(1) = \frac{1}{2} (y+z+1)$ and $f(f(n)) = yn+z$ for $n > 1$.
It is given by
\beql{Eq48}
f\left( c_1 y^k - \frac{z}{y-1} +j \right) = c_2 y^{k+1}
- \frac{z}{y-1} + \frac{y+1}{2} j + \frac{y-1}{2} |j| \,,
\eeq
for $k \ge 0$, where
$$- \frac{y+z-1}{2} y^k \le j < \frac{y+z-1}{2} y^k$$
and
$$c_1 = \frac{(y+1) (y+z-1)}{2(y-1)} , \quad
c_2 = \frac{y+z-1}{y-1} \,.
$$
\end{theorem}

\paragraph{Proof.}
The sequence begins with $f(1) = (y+z+1) /2$, and is constrained by
$f(f(1)) = f((y+z+1)/2) = y+z$.
Since it is monotonically increasing,
the intermediate terms are forced and are given by
$f(1+i) = (y+z+1)/2 +i$ for $0 \le i < (y+z+1)/2$.
The terms $f(2)$, $f(3), \ldots, f((y+z+1)/2)$ now determine the values of
$f(f(2)) = f(y+z)$, $f(f(3)) = f(2y+z) , \ldots$, and we find that the 
stretch from $f((y+z+1)/2)$ through $f(y+z)$ is given by
$f((y+z+1)/2+i) = y+z+yi$ for $0 \le i < (y+z-1)/2$.
Continuing in this way we find that the sequence is completely determined,
and \eqn{Eq48} follows after relabeling the indices.
The conditions \eqn{Eq48a} ensure that there is no 
contradiction in calculating the initial terms, and once 
started the sequence has a unique continuation.~~~$\bsq$

It can be shown (we omit the details) 
that this sequence can also be defined by:
$f(1) = (y+z+1)/2$, and, for $n > 1$, $f(n)$ is the smallest integer
$>f (n-1)$ which is consistent with the condition 
that ``$n$ is a member of the sequence 
if and only if $f(n)$ belongs to the set 
\beql{Eq411}
[ 2, \ldots, \frac{1}{2} (y+z-1) ]
~ \cup ~ \{iy +z : i \ge 1 \} \, . \mbox{''}
\eeq
If $(y+z-1)/2 \le 1$, the first set in \eqn{Eq411} is to be omitted.
\paragraph{Examples.}
Setting $y=3$, $z=0$ in the theorem produces $\{e(n)\}$.

Setting $y=2$, $z=1$ yields another interesting ``mod~2'' sequence.
This is the sequence $\{g(n): n \ge 1 \}$ that begins
$$2,3,5,6,7,9,11,12,13,14, \ldots$$
(\htmladdnormallink{A80637}{http://www.research.att.com/cgi-bin/access.cgi/as/njas/sequences/eisA.cgi?Anum=A080637}).
It has the following properties:

(i)~By definition, this is the unique monotonically
increasing sequence $\{g(n)\}$ satisfying $g(1)=2$,
$g(g(n)) = 2n+1$ for $n \ge 2$.

(ii)~$n$ is in the sequence if and only if $g(n)$ is an odd number $\ge 3$.

(iii)~The sequence of first differences is (\htmladdnormallink{A79882}{http://www.research.att.com/cgi-bin/access.cgi/as/njas/sequences/eisA.cgi?Anum=A079882}):
$$
1, 2, 1^2, 2^2, 1^4, 2^4, 1^8, 2^8, 1^{16}, 2^{16}, \ldots .
$$

(iv)
$$g(3\cdot 2^k -1+j) = 2 \cdot 2^{k+1} -1 + \frac{3}{2} j+ \frac{1}{2} |j| \,,
$$
for $k \ge 0$, $-2^k \le j < 2^k$ (from \eqn{Eq48}).

(v)~$g(2n) = g(n)+g(n-1) +1$, $g(2n+1) = 2g(n) +1$,
for $n \ge 1$ (taking $g(0) =0$).

(vi)~The original sequence $\{a(n)\}$ satisfies
$a(3n) = 3g(n), a(3n+1) = 2g(n)+g(n+1), a(3n+2) = g(n)+2g(n+1)$,
for $n \ge 1$.

(vii)~The ``lying version'' of Section~\ref{Sec3}
is given by $d(n) = g(n+1) -1$ for $n \ge 1$.

(viii)~Let $g'(n) = g(n)+1$.
The sequence $\{g' (n): n \ge 2 \}$ 
was apparently first discovered by C.~L. Mallows,
and is sequence M2317 in \cite{EIS} (\htmladdnormallink{A7378}{http://www.research.att.com/cgi-bin/access.cgi/as/njas/sequences/eisA.cgi?Anum=A007378} in \cite{OEIS}).
This is the unique monotonically increasing
sequence satisfying $g' (g'(n)) = 2n$.
An alternative description is:
$g' (n)$ (for $n \ge 2$) is the 
smallest positive integer $>g' (n-1)$ which is consistent 
with the condition that
\begin{eqnarray}\label{Eqg}
&&\mbox{``$n$ is a member of the sequence if and only if
$g'(n)$ is an even number $\ge 4$''}.
\end{eqnarray}
Note that, although \eqn{Eqc} and \eqn{Eqg} are similar, 
the resulting sequences
$\{c(n)\}$ and $\{g' (n)\}$ are quite different.
$g'$ is not directly covered by Theorem \ref{Th1}, and
we admit that we have not been able 
to identify the largest family of sequences
which can be described by formulae
like
\eqn{Eq1}, \eqn{Eq20},
\eqn{Eq47}, \eqn{Eq48}.

The sequence $\{h(n)\}$ defined by: $h(1) =2$, and, 
for $n > 1$, $h(n)$ is the smallest 
positive integer $>h(n-1)$ which is consistent with the 
condition that ``$n$ is a member of the 
sequence if and only if $h(n)$ is a multiple of 6'':
$$2,6,7,8,9,12,18,24,30,31, \ldots \,,$$
(\htmladdnormallink{A80780}{http://www.research.att.com/cgi-bin/access.cgi/as/njas/sequences/eisA.cgi?Anum=A080780}), shows that such simple rules do not hold in general.
We can characterize the sequence of first differences
in a manner similar to \eqn{Eq1s}, \eqn{Eq1t}:
the alphabet is $\{1, 2, \ldots , 6\}$, and
we define morphisms $\theta (i) = 1, 1, \ldots , 1, 7-i$
(with $i-1$ $1$'s followed by $7-i$), for $i=1, \ldots, 6$. Then the sequence 
of differences of $\{h(h)\}$ is $S_0, S_1, S_2, \ldots $,
where $S_0 = \{4\}, S_{k+1} = \theta(S_k)$ for $k \ge 0$.
However, it appears that no formula similar to \eqn{Eq1}
holds for $h(n)$.

We end this ``mod~m'' section with two interesting ``mod~4''
sequences.
The even numbers satisfy $s(s(n)) = 4n$,
and the odd numbers satisfy $s(s(n)) = 4n+3$.
But there are lexicographically earlier sequences with the 
same properties.
The ``fake even numbers''
$\{i (n) : n \ge 0 \}$ are 
defined by the property that $i(n)$ is the 
smallest nonnegative integer $>i(n-1)$ and 
satisfying $i(i(n)) = 4n$ (\htmladdnormallink{A80588}{http://www.research.att.com/cgi-bin/access.cgi/as/njas/sequences/eisA.cgi?Anum=A080588}):
$$0,2,4,5,8,12,13,14,16,17, \ldots$$
We analyze this sequence by describing the sequence 
of first differences, which are
$$2,2,1,3,4,1,1,2,1,1,1,1,4,4,1,3, \ldots$$
After the initial 2, 2, 1, this breaks up into segments of the form
$$3~S_k~2~T_k ~,$$
where $T_k$ is the reversal of
$$1^1~4^2~1^4~4^8~1^{16}~4^{32} \cdots 4^{2^{2k-1}} 1^{2^{2k}}$$
and $S_k$ is the reversal of
$$1^2~
4^1~1^8~4^2~1^{32}~4^{16} \cdots 1^{2^{2k-1}} 4^{2^{2k-2}} \,.
$$
The ``fake odd numbers'',
$i'(n)$, are similarly defined by $i' (i'(n)) = 4n+3$: 
1, 3, 4, 7, 11, 12, 13, 15, 16, 17, $\ldots$ (\htmladdnormallink{A80591}{http://www.research.att.com/cgi-bin/access.cgi/as/njas/sequences/eisA.cgi?Anum=A080591}), and satisfy
$i'(n) = i(n+1) -1$.

\section{The Aronson transform}\label{Sec5}
A far-reaching generalization of both the original sequence 
and the ``mod~m'' extensions of the previous section 
is obtained if we replace ``odd number''
in the definition of $\{a(n)\}$ by
``member of ${\fbe}$'', where ${\fbe}$ is some fixed sequence.

More precisely, let us fix a starting point $n_0$, which
will normally be $0$ or $1$.
Let ${\fbe} = \{\beta (n) : n \ge n_0 \}$
be an infinite monotonically increasing sequence 
of integers $\ge n_0$ with the property
that its complement (the numbers $\ge n_0$ that are not in ${\fbe}$)
is also infinite.
Then the sequence ${\faf} = \{\alpha (n) : n \ge n_0 \}$
given by: $\alpha (n)$ is
the smallest positive integer $> \alpha (n-1)$ which 
is consistent with the condition that
$$\mbox{``n is in $\faf$ if and only if $\alpha (n)$ is in $\fbe$''}$$
is called the {\em Aronson transform} of ${\fbe}$.

\begin{theorem}\label{ThAT1}
The Aronson transform exists and is unique.
\end{theorem}

\paragraph{Proof.}
For ease of discussion let us call the numbers in  ${\fbe}$ ``hot'',
and those in its complement ``cold.''  We will specify the transform
${\faf}$, leaving to the reader the
easy verification that this has the desired properties,
in particular that there are no contradictions.

The proof is by induction. First we consider the initial term $\alpha(n_0)$.
If $n_0$ is hot, $\alpha(n_0) = n_0$.
If $n_0$ is cold, $\alpha(n_0) = \mbox{~smallest cold number~} \ge n_0 +1$.

For the induction step, suppose $\alpha(n) = k$ for $n > n_0$.

Case (i), $k=n$. 
If $n+1$ is hot then $\alpha(n+1) = n+1$.
If $n+1$ is cold then $\alpha(n+1) = \mbox{~smallest cold number~} \ge n+2$.

Case (ii), $k>n$. 
If $k=n+1$ then $\alpha(n+1) = \mbox{~smallest hot number~} \ge n+2$.
If $k>n+1$ then
if $k+1$ is hot, $\alpha(n+1) = \mbox{~smallest hot number~} \ge k+1$,
while
if $n+1$ is cold,   $\alpha(n+1) = \mbox{~smallest cold number~} \ge k+1$. ~~~$\bsq$

In certain cases it may be appropriate to 
specify some initial terms in $\faf$
to get it started properly.

\paragraph{Examples.}
Of course taking $\fbe$ to be the odd 
numbers (with $n_0=1$) leads to our original sequence $\{a(n)\}$, and 
the even numbers (with $n_0=0$) lead to $\{c(n)\}$ of Section \ref{Sec3}.

The sequences $\PP$ (with $n_0=1$) and $\NN$ (with $n_0=0$)
are fixed under the transformation.

If we take $\fbe$ to be the triangular numbers we get
1, 4, 5, 6, 10, 15, 16, 17, 18, 21, $\ldots$ (\htmladdnormallink{A79257}{http://www.research.att.com/cgi-bin/access.cgi/as/njas/sequences/eisA.cgi?Anum=A079257});
the squares give 1, 3, 4, 9, 10, 11, 12, 13, 16, 25, $\ldots$ (\htmladdnormallink{A79258}{http://www.research.att.com/cgi-bin/access.cgi/as/njas/sequences/eisA.cgi?Anum=A079258});
the primes give 4, 6, 8, 11, 12, 13, 14, 17, 18, 20, $\ldots$ (\htmladdnormallink{A79254}{http://www.research.att.com/cgi-bin/access.cgi/as/njas/sequences/eisA.cgi?Anum=A079254});
and the lower Wythoff sequence (\htmladdnormallink{A201}{http://www.research.att.com/cgi-bin/access.cgi/as/njas/sequences/eisA.cgi?Anum=A000201}),
in which the
$n^{\rm th}$
term is $\lfloor n \phi \rfloor$,
gives 1, 5, 7, 10, 11, 13, 14, 15, 18, 19, $\ldots$ (\htmladdnormallink{A80760}{http://www.research.att.com/cgi-bin/access.cgi/as/njas/sequences/eisA.cgi?Anum=A080760}).

Taking the Aronson transform of $\{a(n)\}$ itself we get
1, 3, 4, 6, 10, 11, 12, 14, 22, 23, $\ldots$ (\htmladdnormallink{A79325}{http://www.research.att.com/cgi-bin/access.cgi/as/njas/sequences/eisA.cgi?Anum=A079325}).
\cite{OEIS} contains several other examples.

The inverse transform may be defined in a similar way.
Given an infinite monotonically increasing sequence 
${\faf} = \{\alpha (n): n \ge n_0 \}$ of numbers $\ge n_0$,
such that its complement (the numbers $\ge n_0$ that are not in ${\faf}$)
is also infinite,
its {\em inverse Aronson transform} is
the sequence ${\fbe} = \{ \beta(n) : n \ge n_0 \}$ 
such that the Aronson transform of ${\fbe}$ is $\faf$.

\begin{theorem}\label{ThAT2}
The inverse Aronson transform exists and is unique.
\end{theorem}

\paragraph{Proof.}
We establish this by 
giving a simple algorithm to construct the inverse transform.
We illustrate the algorithm in Table \ref{Tinv} by 
applying it to the sequence of squares,
${\faf} = \{n^2: n \ge 0 \}$.

Form a table with four rows.
In the first row place the numbers
$n = n_0$, $n_0 +1$, $n_0 +2 , \ldots$, and in 
the second row place the sequence
$\alpha (n_0)$, $\alpha (n_0+1) , \ldots$.
The third row contains what we will call the
``hot'' numbers:
these will comprise the elements of the inverse transform.
The fourth row are the ``cold'' numbers, which 
are the complement of the hot numbers.

The third and fourth rows are filled in as follows.
If $n$ is {\em in} (resp. {\em not in}) the 
sequence $\faf$, place $\alpha (n)$ in the $n$-th
slot of the hot (resp. cold) row.

To complete the table we must fill in the empty slots.
Suppose we are at column $n$, where we 
have placed $\alpha (n)$ in one of the two slots.
Let $l_n$ be the largest number mentioned in 
columns $n_0 , \ldots, n-1$ in the hot or cold rows.
Then we place the numbers $l_n+1 , \ldots 
, \alpha (n) -1$ in the empty slot 
in column $n$, with the single exception that 
if $n$ is not in the sequence and 
$l_n = n-1$ then we place $n$ in the cold slot rather than the hot slot.
(This is illustrated by the position of 2 in the fourth row
of Table \ref{Tinv}.)
We leave it to the reader to verify 
that the entries in the ``hot'' row form 
the inverse Aronson transform $\fbe$. ~~~$\bsq$

\begin{table}[htb]
\caption{Computation of inverse Aronson transform of the squares.  
The ``hot'' numbers comprise the transform.}
$$
\begin{array}{|c|c|c|c|c|c|c|c|c|c|} \hline
n & 1 & 2 & 3 & 4 & 5 & 6 & 7 & 8 & 9 \\ 
\alpha_n & 1 & 4 & 9 & 16 & 25 & 26 & 49 & 64 & 81 \\ 
\mbox{``hot''} & 1 & 3 & \mbox{5--8} & 16 & \mbox{17--24} & \mbox{26--35} & \mbox{37--48} & \mbox{50--63} & 81 \\
\mbox{``cold''} & - & 2,4 & 9 & \mbox{10--15} & 25 & 36 & 49 & 64 & \mbox{65--80} \\ \hline
\end{array}
$$
\label{Tinv}
\end{table}


It follows from Lemma \ref{L1} that $\fbe$ contains 
the members of $\faf^{(2)}$,
but in general $\fbe \neq \faf^{(2)}$.
The additional terms in $\fbe$ serve to make it possible to recover
$\faf$ uniquely from $\fbe$.

\paragraph{Examples.}
As shown in Table \ref{Tinv}, the inverse Aronson transform of the squares (\htmladdnormallink{A10906}{http://www.research.att.com/cgi-bin/access.cgi/as/njas/sequences/eisA.cgi?Anum=A010906}) is
$$1;3;5,6,7,8;
16,17,\ldots,
24; 26, 27, \ldots\,.$$
This consists of a number
of segments (separated here by semicolons).
For $k \ge 1$ the $k^{\rm th}$ segment is $\{k^2\}$ if $k$ is a square,
or $\{ (k-1)^2 +1, \ldots, k^2 -1 \}$ if $k$ is not a square,
except that the second segment $= \{3\}$.

The inverse transform of the primes is
3, 5, 6, 11, 12, 17, 18, 20, 21, 22, $\ldots$ (\htmladdnormallink{A80759}{http://www.research.att.com/cgi-bin/access.cgi/as/njas/sequences/eisA.cgi?Anum=A080759}) 
--- this has a similar decomposition into segments.

The inverse transform of the lower Wythoff sequence is
1, 4, 6, 7, 9, 10, 12, 14, 15, 17, $\ldots$ (\htmladdnormallink{A80746}{http://www.research.att.com/cgi-bin/access.cgi/as/njas/sequences/eisA.cgi?Anum=A080746}).
This consists of the numbers 
$\lfloor \phi k \rfloor + k - 1 ~(k \ge 1)$ and
$\lfloor 2 \phi k \rfloor + k - 1 ~(k \ge 2)$.
The inverse transform of our original sequence $\{a(n)\}$ 
is the sequence of odd numbers (whereas, as we
saw in Section \ref{Sec2}, $a^{(2)}$ omits 3).

In general (because of the above algorithm),
the inverse Aronson transforms are easier to describe 
than the direct transforms.

\section{More complicated conditions}\label{Sec6}
Finally, we may make the condition for $n$ 
to be in the sequence depend on the 
values of several consecutive terms
$a(n)$, $a(n+1), \ldots, a(n+ \tau )$, for some fixed $\tau$.
To pursue this further would take us into the realm of one-dimensional
cellular automata (cf. \cite{Ila01}, \cite{Wol02}), and we will mention
just two examples.

$q(n)$ is the smallest positive integer $> q(n-1)$ 
which is consistent with the condition that ``$n$ 
is in the sequence if and only if 
$q(n)$ is odd and $q(n-1)$ is even'' (\htmladdnormallink{A79255}{http://www.research.att.com/cgi-bin/access.cgi/as/njas/sequences/eisA.cgi?Anum=A079255}):
$$1, 4, 6, 9, 12, 15, 18, 20, 23, 26, 28, \ldots$$
The gaps between successive terms are always 2 or 3.
Changing the condition to ``$\ldots$ both $q(n)$ and $q(n+1)$
are odd'' gives \htmladdnormallink{A79259}{http://www.research.att.com/cgi-bin/access.cgi/as/njas/sequences/eisA.cgi?Anum=A079259}:
$$1, 5, 6, 10, 11, 15, 19, 20, 24, 25,  \ldots$$

\section*{Acknowledgements}
We thank J. C. Lagarias for some helpful comments.


\end{document}